\documentclass[leqno,draft,11pt,a4paper]{article}
\usepackage{amssymb,amsmath,theorem}
\usepackage[includeheadfoot,margin=1in]{geometry}

\allowdisplaybreaks[2]

\newcommand\model{\vDash}
\newcommand\fii{\varphi}
\newcommand\p[1]{\langle#1\rangle}
\DeclareMathOperator\Th{Th}
\DeclareMathOperator\rk{rk}
\newcommand\wo{\mathcal{WO}}
\newcommand\lo{\mathsf{LO}}
\newcommand\ti{\mathsf{TI}}
\newcommand\lango{\mathcal L_<}
\newcommand\EF{\mathrm{EF}}

\begingroup \lccode`\~=`\/ \lowercase{ \endgroup
  \providecommand\url{\begingroup \catcode`\~=12 \catcode`\_=12 \catcode`\/=13 \let~\beginslash \finishurl}
  \newcommand\finishurl[1]{\texttt{#1}\endgroup}
  \newcommand\beginslash{/\futurelet\nexttoken\finishslash}
  \newcommand\finishslash{\ifx\nexttoken~\else\penalty\relpenalty\fi}
}

\newcommand\bme{\hskip.75em\relax}
\newcommand\noproof{\leavevmode\unskip\bme\vadjust{}\nobreak\hfill$\qed$\par}
\newcommand\qed{\Box}
\newenvironment{Pf}[1][]
  {\par\noindent\textit{Proof\optpar{#1}:}\bme\ignorespaces}
  {\noproof\pagebreak[2]\vskip\medskipamount\ignorespacesafterend}
\newcommand\optpar[1]{\ifx\relax#1\relax\else\ #1\fi}

\theoremstyle{plain}
\newtheorem{Thm}{Theorem}
\newtheorem{Lem}[Thm]{Lemma}
\newtheorem{Cl}{Claim}[Thm]

\newenvironment{Pf*}{\let\qed\qedCl\Pf}\endPf

\newcommand\txto{${}\to{}$}

\author{Emil Je\v r\'abek\\[\medskipamount]
Institute of Mathematics, Czech Academy of Sciences\\
\small \v Zitn\'a 25,
115\:67 Praha 1,
Czech Republic,
email: \texttt{jerabek@math.cas.cz}
}

\title{A note on the theory of well orders}

\begin{document}
\maketitle

\begin{abstract}
We give a simple proof that the first-order theory of well orders is axiomatized by transfinite induction, and that it
is decidable.
\end{abstract}

\noindent The first-order theory of the class $\wo$ of well-ordered sets $\p{L,<}$ was developed by Tarski and
Mostowski, and an in-depth analysis was finally published by Doner, Mostowski, and Tarski~\cite{don-mos-tar}: among
other results, they provided an explicit axiomatization for the theory, and proved it decidable. Their main technical
tool is a syntactic elimination of quantifiers, which however takes some work to establish, as various somewhat
nontrivial properties of Cantor normal forms are definable in the theory after all. Alternatively, by way of hammering
nails with a nuke, the decidability of $\Th(\wo)$ follows from an interpretation of the MSO theory of countable linear
orders in the MSO theory of two successors (S2S), which is decidable by a well-known difficult result of
Rabin~\cite{rab:sks}. Our goal is to point out that basic properties of $\Th(\wo)$ can be proved easily using ideas
from L\"auchli and Leonard's~\cite{lau-len} proof of the decidability of the theory of linear orders. A similar
technique was also used by Shelah~\cite{she:mso-lo}.

Let $\lango$ denote the set of sentences in the language $\{<\}$. The theory of (strict) linear orders is denoted
$\lo$; the $\lango$-theory $\ti$ extends $\lo$ with the transfinite induction schema
\[\forall x\:\bigl(\forall y\:(y<x\to\fii(y))\to\fii(x)\bigr)\to\forall x\:\fii(x)\]
for all formulas $\fii$ (that may in principle include other free variables as parameters, though the parameter-free
version is sufficient for our purposes). We will generally denote a linearly ordered set $\p{L,<}$ as just~$L$. Given
linearly ordered sets $I$ and $L_i$ for $i\in I$, let $\sum_{i\in I}L_i$ denote the ordered sum with domain
$\{\p{i,x}:i\in I,x\in L_i\}$ and lexicographic order
\[\p{i,x}<\p{j,y}\iff i<j\text{ or }(i=j\text{ and }x<y).\]
Put $L\cdot I=\sum_{i\in I}L$. We write $L\equiv_kL'$ if $L\model\fii\iff L'\model\fii$ for all $\fii\in\lango$ of
quantifier rank $\rk(\fii)\le k$. It follows from the basic theory of Ehrenfeucht--Fra\"\i ss\'e games (see
\cite{hodges:mod-th,lau-len}) that for each~$k<\omega$, $\equiv_k$ has only finitely many equivalence classes as there
are only finitely many formulas of rank $\le k$ up to logical equivalence, and that $\equiv_k$ preserves sums
and products:
\begin{Lem}\label{lem:sum}
If $L_i\equiv_kL'_i$ for each $i\in I$, then $\sum_{i\in I}L_i\equiv_k\sum_{i\in I}L'_i$. If 
$I\equiv_kI'$, then $L\cdot I\equiv_kL\cdot I'$.
\end{Lem}
\begin{Pf}[sketch]
Duplicator can win the game $\EF_k\bigl(\sum_{i\in I}L_i,\sum_{i\in I}L'_i\bigr)$ by playing auxiliary games
$\EF_k(L_i,L'_i)$ for each $i\in I$ on the side. When Spoiler plays an element in the $i$th summand, Duplicator
simulates it in $\EF_k(L_i,L'_i)$, finds a response using a fixed winning strategy (given by the assumption
$L_i\equiv_kL'_i$), and plays the corresponding element in the main game. The second assertion of the lemma is similar.
\end{Pf}

We come to the main theorem. It was originally proved in~\cite[Thm.~31, Cor.~30, 32]{don-mos-tar} by tedious quantifier
elimination. (The equivalence of \ref{item:1} and~\ref{item:2} also follows from Ehrenfeucht~\cite{ehr:ord}, who proved
by induction on~$k$ that ordinals congruent modulo~$\omega^k$ are $\equiv_k$-equivalent.) Instead, we give a short
argument inspired by the proof of~\cite[Thm.~2]{lau-len} that needs almost no ordinal arithmetic and no explicit EF
game strategies.
\begin{Thm}\label{thm:main}
The following are equivalent for all $\fii\in\lango$:
\begin{enumerate}
\item\label{item:1} $\wo\model\fii$.
\item\label{item:2} $\alpha\model\fii$ for all $\alpha<\omega^\omega$.
\item\label{item:3} $\ti\vdash\fii$.
\end{enumerate}
\end{Thm}
\begin{Pf}
Clearly, \ref{item:3}\txto\ref{item:1}\txto\ref{item:2}. For \ref{item:2}\txto\ref{item:3}, if $\ti\nvdash\fii$,
let $L$ be a countable model of $\ti+\neg\fii$, and $k=\rk(\fii)$; it suffices to show that
there exists $\alpha<\omega^\omega$ such that $L\equiv_k\alpha$. Put
\[S=\left\{c\in L:\forall a,b\in L\:\bigl(a<b\le c\to\exists\alpha<\omega^\omega\:[a,b)\equiv_k\alpha\bigr)\right\}.\]
While the definition speaks of half-open intervals $[a,b)$, the conclusion also holds for $[a,b]$: if
$[a,b)\equiv_k\alpha$, then $[a,b]\equiv_k\alpha+1$ by Lemma~\ref{lem:sum}. Clearly, $S$ is an initial segment of~$L$, and
$0\in S$, where $0=\min(L)$ (which exists by $L\model\ti$).

\begin{Cl}\label{cl:def}
$S$ is definable in~$L$.
\end{Cl}
\begin{Pf*}
Since there are only finitely many formulas of rank $\le k$ up to equivalence, we can form
$\theta_k=\bigwedge\{\theta\in\lango:\rk(\theta)\le k,\forall\alpha<\omega^\omega\,\alpha\model\theta\}$. Then for any linear
order~$L'$, $L'\equiv_k\alpha$ for some $\alpha<\omega^\omega$ iff $L'\model\theta_k$. In particular, $c\in S$ iff
$L\model\forall x,y\,\bigl(x<y\le c\to\theta_k^{[x,y)}\bigr)$, where $\theta_k^{[x,y)}$ denotes $\theta_k$ with quantifiers
relativized to~$[x,y)$.
\end{Pf*}

First, assume that $S$ has a largest element, say~$m$. If $S=L$, then $L=[0,m]\equiv_k\alpha$ for some
$\alpha<\omega^\omega$, and we are done. Otherwise, we will derive a contradiction by showing that the successor of~$m$
(which exists by~$\ti$), denoted~$c$, belongs to~$S$. Indeed, if $a<b\le c$, then either $b\le m$, in which case
$[a,b)\equiv_k\alpha$ for some $\alpha<\omega^\omega$ as $m\in S$, or $b=c$, in which case $[a,b)=[a,m]\equiv_k\alpha$
for some $\alpha<\omega^\omega$ as well.

Thus, we may assume that $S$ has no largest element. Put $S_{\ge a}=\{b\in S:b\ge a\}$.
\begin{Cl}\label{cl:ram}
For every $a\in S$, there is $\alpha<\omega^\omega$ such that $S_{\ge a}\equiv_k\alpha$.
\end{Cl}
\begin{Pf*}
We use the idea of \cite[Lem.~8]{lau-len}. Let $a<a_0<a_1<a_2<\cdots$ be a cofinal sequence in~$S$. For
each $n<m<\omega$, let $t(\{n,m\})=\min\{\alpha<\omega^\omega:[a_n,a_m)\equiv_k\alpha\}$. Since $\equiv_k$ has only
finitely many equivalence classes, $t$ is a colouring of pairs of natural numbers by finitely many colours; by Ramsey's
theorem, there is $\beta<\omega^\omega$ and an infinite $H\subseteq\omega$ such that $t(\{n,m\})=\beta$ for all
$n,m\in H$, $n\ne m$. Let $\{b_n:n<\omega\}$ be the increasing enumeration of $\{a_n:n\in H\}$, and
$\alpha<\omega^\omega$ be such that $[a,b_0)\equiv_k\alpha$. Then
$S_{\ge a}=[a,b_0)+\sum_{n<\omega}[b_n,b_{n+1})\equiv_k\alpha+\beta\cdot\omega<\omega^\omega$
by Lemma~\ref{lem:sum}.
\end{Pf*}

Now, if $S=L$, then $L=S_{\ge0}\equiv_k\alpha$ for some $\alpha<\omega^\omega$ by Claim~\ref{cl:ram}. Otherwise, there
exists $c=\min(L\smallsetminus S)$ by Claim~\ref{cl:def} and $L\model\ti$. We again derive a contradiction by showing
$c\in S$: if $a<b\le c$, then either $b<c$ and $[a,b)\equiv_k\alpha$ for some $\alpha<\omega^\omega$ as $b\in S$, or
$b=c$ and $[a,b)=S_{\ge a}\equiv_k\alpha$ for some $\alpha<\omega^\omega$ by Claim~\ref{cl:ram}.
\end{Pf}

We have so far not actually used any results of L\"auchli and Leonard~\cite{lau-len}, only their methods. But we will
do so now: in order to prove the decidability of $\Th(\wo)$, we need
\begin{Lem}\label{lem:sat}
The relation $\{\p{\alpha,\fii}\in\omega^\omega\times\lango:\alpha\model\fii\}$ is recursively enumerable (thus decidable).
\end{Lem}
Here, we assume $\alpha<\omega^\omega$ is represented by a finite string describing its Cantor normal form (CNF) in a
natural way. Lemma~\ref{lem:sat} is a special case of \cite[Thm.~1]{lau-len}: more generally, L\"auchli and Leonard
prove uniform decidability of linear order types described by ``terms'' using a constant $1$, a binary function~$+$,
unary functions $x\cdot\omega$ and $x\cdot\omega^*$, and a certain variable-arity shuffle operation. It is easy to see
that the CNF of an ordinal $\alpha<\omega^\omega$ can be transformed to such a term using $1$, $+$, and $x\cdot\omega$,
hence Lemma~\ref{lem:sat} follows.

We include a proof of Lemma~\ref{lem:sat} to make the paper more self-contained. It turns out that for well orders,
it is more convenient to consider terms using $1$, $+$, and $\omega\cdot x$ rather than $x\cdot\omega$: then we can
directly axiomatize the theory by a finite sentence without expanding the language with extra predicates as
in~\cite{lau-len}. This argument can be found e.g.\ in Rosenstein~\cite{rost:lo}.

\begin{Pf}[of Lemma~\ref{lem:sat}]
Given $\alpha<\omega^\omega$, we can compute an $\lango$-sentence $T_\alpha$ such that $\Th(\alpha,<)=\lo+T_\alpha$ by
induction on~$\alpha$ as follows. We can take $\forall x,y\,x=y$ for~$T_1$. It is well known that $T_\omega$ can be
defined by axioms postulating that a least element exists, every element has a successor, and every non-minimal element
has a predecessor. If $T_\alpha$ and~$T_\beta$ are already constructed, let $T_{\alpha+\beta}$ be
$\exists x\,\bigl(T_\alpha^{<x}\land T_\beta^{\ge x}\bigr)$, where $T_\alpha^{<x}$ denotes the sentence $T_\alpha$ with all
quantifiers relativized to $(-\infty,x)=\{y:y<x\}$, and similarly for $T_\beta^{\ge x}$: in particular, any
$L\model\lo+T_{\alpha+\beta}$ contains an element $a\in L$ such that $(-\infty,a)\equiv\alpha$ and
$[a,\infty)\equiv\beta$, which implies $L\equiv\alpha+\beta$ by Lemma~\ref{lem:sum}.

Finally, we consider $\omega\cdot\alpha$ for a limit~$\alpha$. Let $\lambda(x)$ denote the formula
$\forall y<x\,\exists z<x\,y<x$, meaning ``$x$ is not a successor''. We define $T_{\omega\cdot\alpha}$ as the
conjuction of $T_\alpha^\lambda$ and axioms postulating that for each $x$, $\max\{y\le x:\lambda(y)\}$ and
$\min\{y>x:\lambda(y)\}$ exist. Clearly, $\omega\cdot\alpha\model T_{\omega\cdot\alpha}$. Conversely, if
$L\model\lo+T_{\omega\cdot\alpha}$, then $L^\lambda:=\{x\in L:L\model\lambda(x)\}\model T_\alpha$, thus
$L^\lambda\equiv\alpha$, and $L=\sum_{x\in L^\lambda}[x,x^+)$, where $x^+=\min\{y>x:y\in L^\lambda\}$. It is easy to
see that $[x,x^+)\model T_\omega$ for each~$x$, thus
$L\equiv\sum_{x\in L^\lambda}\omega=\omega\cdot L^\lambda\equiv\omega\cdot\alpha$ by Lemma~\ref{lem:sum}.
%
\end{Pf}

The following consequence is Theorem~33 of~\cite{don-mos-tar}.
\begin{Thm}\label{thm:dec}
The theory $\Th(\wo)=\ti$ is decidable.
\end{Thm}
\begin{Pf}
$\{\fii\in\lango:\ti\vdash\fii\}$ is recursively enumerable as $\ti$ is recursively axiomatized; by
Theorem~\ref{thm:main} and Lemma~\ref{lem:sat},
$\{\fii\in\lango:\ti\nvdash\fii\}=\{\fii\in\lango:\exists\alpha<\omega^\omega\,\alpha\model\neg\fii\}$ is also
recursively enumerable.
\end{Pf}

\paragraph{Acknowledgements.}
I am grateful to the anonymous reviewers for helpful suggestions, especially the simplification of the proof of
Lemma~\ref{lem:sat}.

The author was supported by the Czech Academy of Sciences (RVO 67985840) and GA \v CR project 23-04825S. The
main argument was developed while teaching a course on ``Decidable theories'' at the Charles University in Prague
(\url{https://users.math.cas.cz/\string~jerabek/teaching/decidable.html}).

\bibliographystyle{mybib}
\bibliography{wo}

\end{document}